\documentclass[11pt]{amsart}
\usepackage{latexsym}  
\usepackage{amssymb,amscd,amsthm,amsfonts,verbatim,amsopn,color}
\usepackage[dvips]{epsfig}
\def\Q{{\mathbb Q}}
\def\C{{\mathbb C}}

\def\cd{{\rm cd}}
\def\deg{{\rm deg\,}}

\newcommand{\LL}{\mathcal L}

\newcommand{\CC}{\mathcal C}

\newcommand{\Aa}{\mathcal A}

\newcommand{\MM}{\mathcal M}

\newcommand{\rk}{{\rm rk}}

\newtheorem{thm}{Theorem}[section]
\newtheorem{df} [thm]{Definition}

\newtheorem{prop}[thm]{Proposition}

\newtheorem{rem}[thm]{Remark}

\newtheorem{expl}[thm]{Example}
\numberwithin{equation}{section}

\newenvironment{pf}{\noindent {\bf Proof.}}{\hfill $\Box$\vspace{0.3cm}}
\newenvironment{explrm}{\begin{expl} \rm}{\end{expl}}
\newenvironment{remrm}{\begin{rem} \rm}{\end{rem}}
\newenvironment{dfrm}{\begin{df} \rm}{\end{df}}

\begin{document}

\title[Formality of geometric arrangements]
{Formality of the complements of subspace arrangements 
with geometric lattices}

\author{Eva Maria Feichtner \mbox{ }and \mbox{ }Sergey Yuzvinsky}

\address{
Department of Mathematics, ETH Zurich, 8092 Zurich, Switzerland}
\email{feichtne@math.ethz.ch}

\address{
Department of Mathematics, University of Oregon, Eugene OR 97403-1222, 
USA}
\email{yuz@math.uoregon.edu}


\begin{abstract} 
  We show that, for an arrangement of subspaces in a complex vector
  space with geometric intersection lattice, the complement of the
  arrangement is formal. We prove that the Morgan rational model for
  such an arrangement complement is formal as a differential graded
  algebra.
\end{abstract}


\maketitle


\section{Introduction} \label{sect_intr}

Let $\Aa$ be an arrangement of linear subspaces in a complex vector 
space~$V$; we denote the complement of $\Aa$ by 
$\MM(\Aa)\,{:=}\,V\,{\setminus}\,\bigcup \Aa$. In this note we address the 
question whether  $\MM(\Aa)$ is a formal space. If $\Aa$ is an arrangement
of hyperplanes then formality of the arrangement complement follows 
immediately from Brieskorn's result \cite[Lemme 5]{B73}, which gives an 
embedding of $H^*(\MM(\Aa))$ into the deRham complex of~$\MM(\Aa)$. 
In the general case, however, the question is more subtle.

Proving formality can be reduced to an algebraic question.  Indeed,
using results of Morgan from~\cite{M78}, to prove that $\MM(\Aa)$ is
formal it suffices to prove that the Morgan rational model of
$\MM(\Aa)$ is formal as a differential graded algebra.

Formality of a differential graded algebra, i.e., the property that it
is quasi-isomorphic to its cohomology algebra, is preserved under the
equivalence relation on differential graded algebras that is generated
by quasi-isomorphisms. Using this fact, we construct a sequence of
quasi-isomorphic differential graded algebras, starting from the
Morgan rational model of $\MM(\Aa)$ and finishing on a differential
graded algebra with a quasi-isomorphism to its cohomology algebra.

Morgan rational models for arrangement complements have not only been
described explicitly in work of De~Concini and Procesi~\cite{DP95},
there are also a number of quasi-isomorphic variations at hand. The
first quasi-isomorphism that we rely on is the quasi-isomorphism
between the De~Concini-Procesi rational model for arrangement
complements and a differential graded algebra $\mathrm{CM}_{\Aa}$
introduced by the second author in~\cite{Yu02}.  Its underlying chain
complex is the flag complex of the intersection lattice of $\Aa$.  We
go further and prove that $\mathrm{CM}_{\Aa}$ in turn is
quasi-isomorphic to a differential graded algebra~$D_{\Aa}$, which as
well has been introduced by the second author~\cite{Yu99}. In contrast
to $\mathrm{CM}_{\Aa}$, the underlying chain complex of $D_{\Aa}$ is
the relative atomic complex of the intersection lattice of $\Aa$.

It is for $D_{\Aa}$ that we prove formality in the case that the
intersection lattice of the arrangement is a geometric lattice.  Only
in this last part of our work, we have to rely on specific properties
of geometric intersection lattices. We remark here that formality of
arrangements with geometric intersection lattices was mentioned without a
proof in \cite[Remark 7.3(ii)]{Yu02}, and was used recently 
in work of Papadima and Suciu~\cite{PS04}.
Our goal here is to present all the details of the proof.

\medskip

The paper is organized as follows: We start out by recalling the
definition of quasi-isomorphism and formality for differential graded
algebras.  We then introduce the central for our approach object, the
differential graded algebra $D_{\Aa}$ for a given subspace
arrangement.  In Section~\ref{ssect_arrgtcoh} we explain its role for
arrangement cohomology by showing that it is quasi-isomorphic to the
arrangement model $\mathrm{CM}_{\Aa}$ constructed in~\cite{Yu02}.
Section~\ref{sect_formal} finally is devoted to proving that $D_{\Aa}$
is formal if the arrangement has a geometric intersection lattice. In
the last section, we discuss the perspectives and limitations of our
approach towards proving formality for more general subspace
arrangements.


\section{Preliminaries}  \label{sect_prelim}

\subsection{Formality of differential graded algebras}
\label{ssect_formality}

We recall the definitions of quasi-isomorphisms of 
differential graded algebras (d.g.a.), and of the equivalence relation
they generate. Henceforth, we give the concept of formality
that we will use in the sequel.

\begin{dfrm} \label{def_qi}
Let $(E,d_E)$ and $(F,d_F)$ be differential graded algebras.
A d.g.a. morphism $h:(E,d_E) \rightarrow 
(F,d_F)$ is called 
{\em quasi-isomorphism\/} if the induced map in cohomology, 
$h^*: H^*(E)\rightarrow H^*(F)$, is an isomorphism.
\end{dfrm}

The existence of a quasi-isomorphism from one differential graded
algebra to another defines a relation which is not in general
symmetric or transitive. We consider the equivalence relation of
d.g.a.'s that is generated by quasi-isomorphisms, and we say that two
d.g.a.'s $E$ and $F$ are {\em quasi-isomorphic\/} if there exists a
finite sequence of d.g.a.'s starting on $E$ and finishing on $F$, each
neighboring pair being connected by a quasi-isomorphism in at least
one direction.

\begin{dfrm} \label{def_formal}
A differential graded algebra $(E,d_E)$ is called {\em formal\/} 
if it is quasi-isomorphic to its cohomology  algebra $H^*(E)$ with 
zero differential.   
\end{dfrm}

Notice that formality of a differential graded algebra 
is invariant under the relation of being quasi-isomorphic.



\subsection{The relative atomic differential graded algebra of an arrangement}
\label{ssect_atomicd.g.a.}

Let~$\Aa$ be an arrangement of linear subspaces in a complex vector
space, and let $\LL$ denote its intersection lattice, i.e., the poset
of intersections among subspaces in $\Aa$ ordered by reversed
inclusion. We will frequently refer to the (complex) codimension, $\cd
\, A$, of elements $A$ in $\LL$ as subspaces in the ambient space of
$\Aa$. We assume that spaces in~$\Aa$ are inclusion maximal, hence in
one-to-one correspondence with the atoms $\mathfrak A(\LL)$ in $\LL$.
Furthermore, we fix a linear order on the set of atoms $\mathfrak
A(\LL)$.

We define the {\em relative atomic differential graded algebra\/}
 $D_{\Aa}$ associated with an arrangement~$\Aa$ as follows. The
underlying chain complex $(D,d)$ is the {\em relative atomic complex\/},
with coefficients in~$\Q$. We recall its definition 
from~\cite[Def.\ 2.2]{Yu99}. The complex $(D,d)$ is generated by all 
subsets $\sigma$ of
$\mathfrak A(\LL)$, and for $\sigma\,{=}\,\{i_1,\ldots, i_k\}$ in
$\mathfrak A(\LL)$, the differential $d$ is defined by
\begin{equation} \label{eq_dsigma}
                     d\sigma\, = \,
                     \sum_{j: \bigvee \sigma_j = \bigvee \sigma}
                    (-1)^j\, \sigma_j                                
\end{equation}  
where $\sigma_j\,{=}\,\sigma\setminus\{i_j\}$ for $j=1,\ldots,k$, and the 
indexing of elements in $\sigma$ follows the linear order imposed on 
$\mathfrak A(\LL)$. With 
$\deg(\sigma)\,{=}\, 2 \,\cd \bigvee \sigma - |\sigma|$, 
$(D,d)$ is a cochain complex.

We define a multiplication on $(D,d)$ as follows. For subsets $\sigma$ and 
$\tau$ in $\mathfrak A(\LL)$,
\begin{equation} \label{eq_multD}
    \sigma \cdot \tau \,\, = \,\,\left\{  \begin{array}{cl}
    (-1)^{{\rm sgn} \epsilon(\sigma,\tau)} \sigma \cup \tau& 
    \mbox{if }\, \cd \bigvee \sigma + \cd \bigvee \tau = 
                              \cd \bigvee (\sigma \cup \tau)  \\
                             0       & \mbox{otherwise}
                               \end{array}
                      \right. \, , 
\end{equation}
where $\epsilon(\sigma,\tau)$ is the permutation that, applied to
$\sigma \cup \tau$ with the induced linear order, places elements of
$\tau$ after elements of $\sigma$, both in the induced linear order.

\begin{thm} \label{thm_d.g.a.D}
{\rm (\cite[Prop.\ 3.1]{Yu00})} For any arrangement $\Aa$ of complex subspaces,
$D_{\Aa}$ with underlying chain complex $(D,d)$  and multiplication defined
in~{\rm (\ref{eq_multD})} is a differential graded algebra.  
\end{thm}

\begin{remrm} \label{rem_Dd} \mbox{ } 
{\bf (1)} A differential graded
algebra similar to the one discussed here, can be defined for arbitrary
lattices with a labelling of elements that satisfies certain rank-like
conditions. For a detailed discussion of the general context, 
see~\cite[Section~3]{Yu00}.
\newline 
{\bf (2)} Recall that there are two abstract simplicial
  complexes associated with any finite lattice $\LL$: the {\em flag
    complex\/} or {\em order complex\/} $F(\LL)$ and the {\em atomic
    complex\/} $C(\LL)$.  The flag complex is formed by all flags,
  i.e., linearly ordered subsets, in the reduced lattice $\bar \LL{=}
  \LL\setminus \{\hat 0, \hat 1\}$.  The atomic complex consists of
  all subsets $\sigma$ of $\mathfrak A(\LL)$ with $\bigvee
  \sigma\,{<}\,\hat 1$.  The abstract simplicial complexes $F(\LL)$
  and $C(\LL)$, in fact, are homotopy equivalent.  To simplify
  notation, we will not distinguish between an abstract simplicial
  complex and its simplicial chain complex with rational coefficients.
  
We can now give an explanation for the terminology chosen for the
  differential graded algebra $D_{\Aa}$. Observe that the complex
  $(D,d)$ naturally decomposes as a direct sum
\begin{equation} \label{eq_decD} 
    D \,\, = \,\, \bigoplus_{A\in \LL} D(A)\, ,
\end{equation}
where $D(A)$, for $A\,{\in}\,\LL$, is generated by all subsets
$\sigma$ in $\mathfrak A(\LL)$ with $\bigvee \sigma\,{=}\,A$.

In fact, for any $A\,{\in}\,\LL$, and $p\,{\geq}\,0$, there is a
natural isomorphism
\begin{equation} \label{eq_hDA}
    H^p(D(A)) \, \, \cong \, \, \widetilde H_{2\cd A-p-2} (C(\hat 0, A)) \, ,
\end{equation}
where $C(\hat 0, A)$ denotes the atomic complex of the interval $[\hat
0, A]$ in $\LL$.

It is an easy observation that $(D(A),d)$, graded by cardinality of
generators $\sigma\,{\subseteq}\,\mathfrak A(\LL)$, is the same as the
relative simplicial chain complex $\Sigma(A)/ C(\hat 0,A)$, where
$\Sigma(A)$ is the simplicial chain complex of the full simplex on the
vertex set $\{X\in \mathfrak A(\LL)\,|\, X\leq A\}$. The isomorphism
stated above is part of the exact homology pair sequence of
$(\Sigma(A),C(\hat 0,A))$; compare~\cite[Sect.\ 3.1.2]{Yu01} for
details. We will later write out a chain map that induces the
isomorphism in~(\ref{eq_hDA}) as part of our proof of
Proposition~\ref{prop_DCM}.  \newline 
{\bf (3)} For an arrangement
with geometric intersection lattice, $H^*(D,d)$ is generated by the
classes $[\sigma]$ of independent subsets $\sigma\subseteq \mathfrak
A(\LL)$. Indeed, for any independent set $\sigma$ in $\LL$, $\sigma$
is a cocycle in $(D,d)$ by definition of the differential. Moreover,
from the description of homology in~(\ref{eq_hDA}), and using the
well-known results of Folkman on lattice homology~\cite{Fo66}, we see
that $H^p(D(A))\,{=}\,0$ unless $p\,{=}\, 2 \cd A-\rk A$. For
$\sigma\,{\subseteq}\,\mathfrak A(\LL)$ to be a generator of $D^{2 \cd
  A-\rk A}(A)$, is is a necessary condition that $\bigvee
\sigma\,{=}\,A$ and $|\sigma|\,=\,\rk A$, hence $\sigma$ has to be
independent. It follows that the classes $[\sigma]$, $\sigma$
independent, generate $H^*(D_{\Aa})$.
\end{remrm}


\subsection{Arrangement cohomology}
\label{ssect_arrgtcoh} 
\mbox{ } \newline We here recall the definition of another
differential graded algebra ${\rm CM}_{\Aa}$ associated with any
complex subspace arrangement. This differential graded algebra is the
main character in~\cite{Yu02}, where it is shown to be a rational
model for the arrangement complement. It is a considerable
simplification of the rational model presented earlier in work of
De~Concini and Procesi~\cite{DP95}, and relates their results in an
elucidating way to the much earlier results of Goresky~\&
MacPherson~\cite{GM88} on the linear structure of arrangement
cohomology.

The underlying chain complex of ${\rm CM}_{\Aa}$ is the $\LL$-graded complex
\[
      {\rm CM}_{\Aa}\, \, = \, \, \bigoplus_{A\in \LL} F(\hat 0, A)\, ,
\]
where $F(\hat 0, A)$, for $A\,{\in}\,\LL$, is the flag complex of the 
open interval $(\hat 0,A)$ of elements in $\LL$ below $A$.

To describe the multiplication we need to fix some notation. 
Given a flag $T\,\in\, F(\hat 0,A)$, denote by $\overline T$ the flag 
extended by the lattice element $A$, and by $\underline  T$ the flag with its 
maximal element removed. For an ordered collection of lattice elements
$\CC\,{=}\,\{C_1,\ldots, C_t\}$,  denote by $\lambda(\CC)$ the chain in 
$\LL$ obtained by taking successive joins,
\[
 \lambda(\CC): \, \, C_1 < C_1\vee C_2 < \ldots < 
 \bigvee_{i=1}^t C_i\, ,    
\]
and set $\lambda(\CC)$ to zero in case there are repetitions occurring
among the joins.
 
For $A,B\,{\in }\,\LL$, let $T_A: A_1<\ldots <A_p$, $T_B: B_1<\ldots<
B_q$, be flags in $F(\hat 0,A)$, $F(\hat 0,B)$, respectively. Let
$\mathfrak S_{p+1,q+1}$ denote the shuffle permutations in the
symmetric group $\mathfrak S_{p+q+2}$, i.e., permutations of $[p+q+2]$
that respect the relative order of the first $p{+}1$ and the relative
order of the last $q{+}1$ elements. Denote by $(\overline
T_A,\overline T_B)^{\pi}$ the result of applying $\pi\in \mathfrak
S_{p+1,q+1}$ to the pair of chains, with elements of $\overline T_A$ in
ascending order preceding elements of $\overline T_B$ in ascending
order, and thereafter applying $\lambda$,
\[
(\overline T_A,\overline T_B)^{\pi} 
\, \, = \, \, 
\lambda \pi(\overline T_A,\overline T_B ) \, .
\]
We are now ready to describe the product on ${\rm CM}_{\Aa}$. 
For $T_A$ and $T_B$
as above, we define
\begin{equation} \label{eq_multCM}
    T_A \cdot T_B \,\, = \,\,\left\{  \begin{array}{cl}
{\displaystyle    
\sum_{\pi\in \mathfrak S_{p+1,q+1}} 
    (-1)^{{\rm sgn}\,\pi}\, \,  \underline{(\overline T_A,\overline T_B)^{\pi}}
}& 
    \, \, \, \mbox{if }\, \cd A + \cd B = 
                              \cd (A \vee B)  \\[0.6cm]
                             0       & \mbox{otherwise}
                               \end{array}
                      \right. \, . 
\end{equation}

As we mentioned in the introduction, $\mathrm{CM}_{\Aa}$ is our link 
between the relative atomic differential graded algebra $D_{\Aa}$ and 
the De~Concini-Procesi rational model for arrangement complements given
in~\cite{DP95}.

\begin{thm}{\rm \cite[Cor.\ 4.7 and 5.3]{Yu02}}
The differential graded algebra ${\rm CM}_{\Aa}$ 
associated with an arrangement 
$\Aa$ is quasi-isomorphic to the De~Concini-Procesi rational model 
for the arrangement complement. 
\end{thm}

To complete the sequence of quasi-isomorphisms, 
we are left to show the following.

\begin{prop} \label{prop_DCM}
  The relative atomic differential graded algebra $D_{\Aa}$ of an
  arrangement~$\Aa$ is quasi-isomorphic to the differential graded
  algebra ${\rm CM}_{\Aa}$.
\end{prop}

\begin{pf}
We define a homomorphism of differential graded algebras 
$h:D_{\Aa} \longrightarrow \mathrm{CM}_{\Aa}$. The homomorphism respects 
the $\LL$-grading of both algebras, hence it will be sufficient to define
$h_A:D(A) \longrightarrow F(\hat 0, A)$ for any $A\,{\in}\,\LL$. 
The ingredients 
are two maps, $g_A: D(A) \longrightarrow C(\hat 0,A)$ and 
$f_A: C(\hat 0,A) \longrightarrow F(\hat 0,A)$, where $C(\hat 0,A)$ 
is the atomic complex of the interval $[\hat 0,A]$ in $\LL$. 
For $\sigma\,{\subseteq}\,\mathfrak A(\LL)$ with $\bigvee \sigma\,{=}\,A$, 
we define   
\begin{equation*} 
                     g_A(\sigma)\, = \,
                     \sum_{j: \bigvee \sigma_j < \bigvee \sigma}
                    (-1)^j\, \sigma_j \, ,                               
\end{equation*} 
where, as in the definition of the differential in $D_{\Aa}$ 
in~(\ref{eq_dsigma}), the indexing of elements in $\sigma$ follows the 
linear order imposed on $\mathfrak A (\LL)$. 

The map $f_A$ is the standard chain homotopy equivalence between the 
atomic complex and the flag complex of a given lattice, which we recall 
from~\cite[Lemma 6.1]{Yu02} for completeness. For 
$\sigma\,{\subseteq}\,\mathfrak A(\LL)$ with $\bigvee \sigma\,{<}\,A$,
we define
\begin{equation*} 
                     f_A(\sigma)\, = \,
                     \sum_{\pi \in \mathfrak S_{|\sigma|}} 
                    (-1)^{\mathrm{sgn}\, \pi} \lambda(\sigma)\, .
\end{equation*} 
We now define for 
$\sigma\,{\subseteq}\,\mathfrak A(\LL)$ with $\bigvee \sigma\,{=}\,A$,
\begin{equation} \label{eq_defh}
        h_A (\sigma)
\, \, = \, \, 
        (-1)^{2 \cd A -|\sigma|}\, f_A\, g_A (\sigma). 
\end{equation}
We claim that $h=\sum_{A\in \LL}h_A$ is a quasi-isomorphism of differential 
graded algebras, and we break our proof into 4 steps.

\noindent
{\bf (1)} {\em $h$ is a homogeneous map\/}. This is obvious from the 
definition of the gradings on $D_{\Aa}$ and $\mathrm{CM}_{\Aa}$. \newline
{\bf (2)} {\em $h_A$ is a map of chain complexes\/}. 
Due to the sign in (\ref{eq_defh}), we need to show that 
$f_A\, g_A (d\sigma) \,{=}\, - d( f_A\,g_A(\sigma))$  for 
$\sigma\,{\subseteq}\,\mathfrak A(\LL)$ with $\bigvee \sigma\,{=}\,A$. 
By definition of $d$ and $g_A$, we have
\[
     (d\,+ \, g_A)(\sigma) \, \, = \, \,
                     \sum_{j\in \sigma}
                    (-1)^j\, \sigma_j \, ,     
\]
which implies that $(d{+}g_A)^2\,{=}\,0$. With $d^2\,{=}\,0$ and 
$g_A^2\,{=}\,0$, we conclude that $dg_A\,{=}\,-g_Ad$, and since $f_A$ 
is a chain map, our claim follows.\newline
{\bf (3)} {\em $h$ is multiplicative\/}.  Let $\sigma$, $\tau$ be subsets 
in $\mathfrak A(\LL)$ and denote $\bigvee \sigma\,{=}\,A$ and  
$\bigvee \tau\,{=}\,B$. We can assume that $\cd A + cd B=\cd (A\vee B)$, 
in particular, $\sigma\cap \tau\,{=}\,\emptyset$, since otherwise both 
sides of the equation  
$h_{A\vee B}(\sigma \cdot \tau)\,{=}\, h_{A}(\sigma)\cdot h_{B}(\tau)$ are $0$ 
by definition of  products in $D_{\Aa}$ and $\mathrm{CM}_{\Aa}$, respectively. 
The sign in the definition of $h$ is chosen so that here we 
are left to show that 
\begin{equation} \label{eq_fgmult}
  f_{A\vee B}g_{A\vee B}(\sigma\cdot \tau)
  \, \, = \, \, 
  f_Ag_A(\sigma)\cdot f_Bg_B(\tau)\, .
\end{equation}
Now we make two claims. \newline {\bf Claim 1.} The set of nonzero flags
occurring on the left hand side of~(\ref{eq_fgmult}) coincides with
the set of flags on the right hand side. \newline 
{\bf Claim 2.} The
signs of a flag occurring on the left hand side and on the right hand
side of~(\ref{eq_fgmult}) are the same.

 To prove Claim 1, fix a nonzero flag $F$ occurring on the left hand
 side of~(\ref{eq_fgmult}). It is constructed by applying $\lambda$  
 to some linear order
 on $(\sigma\cup\tau)\,{\setminus}\,\{a\}$ for an atom $a$.  Suppose
 $a\,{\in}\,\sigma$, the other case being similar. Denote by $b$ the
 last element of $\tau$ in this ordering. Since $F$ does not have
 repetitions, we find $\pm\tau\,{\setminus}\,\{b\}$ as a summand in
 $g_B(\tau)$. The induced orderings on $\sigma\,{\setminus}\,\{a\}$
 and $\tau\,{\setminus}\,\{b\}$ produce flags $F_1$ and $F_2$ in
 $F(\hat 0,A)$ and $F(\hat 0,B)$, respectively, whose product occurs
 on the right hand side of~(\ref{eq_fgmult}). The linear order that
 gave rise to $F$ prescribes a shuffle permutation that generates $\pm
 F$ as a summand of $F_1\cdot F_2$ on the right hand side.  The
 opposite inclusion can be shown by inverting all the steps of the
 proof.
 
 To prove Claim 2, fix again a flag $F$ occurring on the left hand
 side of~(\ref{eq_fgmult}) as above, and keep the notation from the
 proof of Claim 1.  Denote the coefficients of $F$ on the left hand
 side and on the right hand side by $(-1)^{\ell}$ and $(-1)^r$,
 respectively.

Then we have
\begin{equation} \label{eq_claim2_lhs}
\ell\, \, =\mathrm{sgn}\,\epsilon(\sigma,\tau)
\,+\,[a\in\sigma\cup\tau]
\,+\,[(\sigma\setminus\{a\})\cup\tau]\, ,
\end{equation}
where the second term is the numerical position of $a$ in
$\sigma\cup\tau$ (in the initial ordering) and $[\rho]$ is the parity
of the permutation induced by the ordering that gives rise to $F$ on
any subset $\rho$ of $(\sigma\setminus\{a\})\cup\tau$. Using similar
notation, we have
\begin{equation} \label{eq_claim2_rhs}
r\,\, =\,\, [a\in\sigma]
\,+\,[b\in\tau]
\,+\,[\sigma\setminus\{a\}]
\,+\,[\tau\setminus\{b\}]
\,+\,[\sigma\setminus\{a\},\tau]
\,+\,|\tau|\, ,
\end{equation}
where $[\rho_1,\rho_2]$ is the parity of the shuffle permutation, 
induced by the
ordering that gives rise to $F$, of two disjoint subsets $\rho_1$ and
$\rho_2$ of $(\sigma\setminus\{a\})\cup(\tau\setminus\{b\})$ (from the
starting position of $\rho_2$ after $\rho_1$). Recall that, to obtain
$F$, we need first to augment $f_A(\sigma\setminus\{a\})$ by $A$, and
$f_B(\tau\setminus\{b\})$ by $B$, respectively, and then apply the
needed shuffle. It is easy to see that the shuffle should have $A$ at
the end of the set. Thus the last summand in~(\ref{eq_claim2_rhs})
comes from moving $A$ over $\tau$. The augmentation by $B$ amounts
just to the substitution of $b$ by $B$ in $\tau$.

Now we need the following straightforward equalities (modulo 2).
\begin{itemize}
\item[(i)]$\mathrm{sgn}\,\epsilon(\sigma,\tau)\,{=}\,
  \mathrm{sgn}\,\epsilon(\sigma\setminus\{a\},\tau)
  \,{+}\,|\tau_{<a}|\,$, \newline 
  where the new symbols are self-explanatory (e.g.,
  $\tau_{<a}\,{=}\,\{c\in\tau| c<a\}$);
\item[(ii)]
$[a\in\sigma\cup\tau]\,{=}\,[a\in\sigma]\,{+}\,|\tau_{<a}|$;
\item[(iii)]
$[\sigma\setminus\{a\}]\,{+}\,[\tau\setminus\{b\}]\,{+}\,
[\sigma\setminus\{a\},\tau]\,{=}\,
\mathrm{sgn}\,\epsilon(\sigma\setminus\{a\},\tau)\,{+}\,|\tau_{>b}|\,{+}\,
[(\sigma\setminus\{a\})\cup\tau]$;
\item[(iv)] $[b\in\tau]\,{+}\,|\tau_{>b}|\,{=}\,|\tau|$.
\end{itemize}
Substituting (i) - (iv) in~(\ref{eq_claim2_lhs}) and~(\ref{eq_claim2_rhs}) 
we obtain the needed equality of the
coefficients of $F$ on the left and right hand sides of~(\ref{eq_fgmult}).

{\bf (4)}
{\em $h$ is a quasi-isomorphism\/}. 
Both maps $g_A$ and $f_A$ induce isomorphisms in homology, which 
completes our proof.
\end{pf}

To prove that the complement of an arrangement $\Aa$ is formal,
we are left to show that the relative atomic differential graded
algebra $(D_{\Aa},d)$ is formal.  Note that so far (with the exception
of Remark~\ref{rem_Dd}(3)) we did not refer to any specific properties
of the arrangement or its intersection lattice. It is only in the next
section that we will restrict ourselves to arrangements with geometric
intersection lattices.


\section{Formality of $D_{\Aa}$ for geometric lattices}  
\label{sect_formal}

\begin{thm} \label{thm_Dformal}
Let $\Aa$ be a complex subspace arrangement with geometric intersection 
lattice. 
The linear map
\begin{eqnarray*}
\Psi\,:\, \, D_{\Aa} & \longrightarrow & H^*(D_{\Aa},d) \\
\sigma &\longmapsto & \left\{  \begin{array}{cl}
        [\sigma], & \mbox{if }\, \sigma \mbox{ is independent}\, ,\\
                             0,       & \mbox{otherwise}\, .
                               \end{array}
                      \right. 
\end{eqnarray*}
is a quasi-isomorphism of differential graded algebras, i.e.,
$D_{\Aa}$ is formal. 
\end{thm}

\begin{pf}
{\bf (1)} 
{\em $\Psi$ is multiplicative.\/} We need to check that for any two
subsets $\sigma$, $\tau$ in $\mathfrak A(\LL)$, 
\begin{equation} \label{eq_mult_psi}
       \Psi(\sigma \tau) \, = \, \Psi(\sigma) \Psi(\tau)\, . 
\end{equation}
\noindent
(i) First assume that both $\sigma$ and $\tau$ are independent in
$\LL$, and $\cd \bigvee \sigma \,{+}\, \cd \bigvee \tau \,{=}\, \cd
\bigvee \sigma \cup \tau$. \newline 
{\em Claim.\/} $\sigma \cup \tau$ is independent in $\LL$.\newline 
We first observe that $\sigma \cap
\tau =\emptyset$. For, if $c\in \sigma \cap \tau$, then 
$\cd((\bigvee \sigma) \wedge (\bigvee \tau)) \geq \cd \,c > 0$, and 
\[
\cd \bigvee (\sigma\cup \tau) < 
\cd \bigvee (\sigma\cup \tau)+\cd((\bigvee \sigma) \wedge (\bigvee \tau)) \leq
\cd \bigvee \sigma + \cd \bigvee \tau\,,
\]
contrary to our assumption. 

Assume that $\sigma{\cup}\tau$ were dependent, hence there existed a
$c\,{\in}\,\sigma \,{\cup}\,\tau$ with $\bigvee (\sigma{\cup}\tau
\,{\setminus}\,\{c\}) = \bigvee (\sigma{\cup}\tau)$. Without loss of
generality, we can assume that $c\in \tau$. Then
\[
\cd \bigvee (\sigma\cup\tau) = 
\cd ( \bigvee \sigma \,\vee \,\bigvee(\tau \setminus \{c\})) \leq 
\cd \bigvee \sigma + \cd \bigvee (\tau \setminus \{c\}) < 
\cd \bigvee \sigma + \cd \bigvee \tau\, ,
\]
 where the strict inequality holds since
$\tau$ is independent. Again, we reach a contradiction to our assumption.

With $\sigma\,{\cup}\,\tau$ being independent, equality~(\ref{eq_mult_psi})
holds by definition of multiplication in $H^*(D_{\Aa},d)$.

\noindent
(ii) Now assume that $\sigma$ and $\tau$ are independent in
$\LL$, but $\cd \bigvee \sigma \,{+}\, \cd \bigvee \tau \,{\neq}\, \cd
\bigvee (\sigma \cup \tau)$. \newline
In this case $\sigma \tau =0$ in $D$, hence both sides in~(\ref{eq_mult_psi})
are zero, the right hand side again by definition of multiplication 
in $H^*(D,d)$.

\noindent
(iii) We conclude by assuming that at least one of $\sigma$ or $\tau$
are dependent sets, say $\sigma$. With $\psi(\sigma)=0$ by definition,
the right hand side in~(\ref{eq_mult_psi}) is zero. As $\sigma$ is
dependent, so is $\sigma\cup \tau$. Either $\sigma\tau=0$ in $D$, then
so is its image under $\Psi$, or codimensions add up and $\sigma\tau$
equals $\sigma\cup \tau$ up to sign. Then, the left hand side is zero
by definition of~$\Psi$.

\smallskip
\noindent
{\bf (2)} 
{\em $\Psi$ is a homomorphism of differential graded algebras.\/}
We need to show that $\Psi(d\sigma)=0$ for any subset $\sigma$ in 
$\mathfrak A(\LL)$.
\newline
For an independent set $\sigma$ in $\LL$, $d\sigma=0$ by definition, hence 
$\Psi(d\sigma)=0$. \newline
Let $\sigma$ be a dependent set in $\LL$. We recall that 
\begin{equation*} 
                     d\sigma\, = \,
                     \sum_{j: \bigvee \sigma_j = \bigvee \sigma}
                    (-1)^j\, \sigma_j \, ,                               
\end{equation*}     
and we observe that $d\sigma$ maps to zero under $\Psi$ if either all
summands $\sigma_j$ in (\ref{eq_dsigma}) are dependent or if all
summands in (\ref{eq_dsigma}) are independent.  In the latter case,
$\Psi(d\sigma)=[d\sigma]$ is the homology class of a boundary in $D$,
hence is zero. We need to prove that these two cases exhaust all
possibilities.

Assume that $\sigma_i$ is an independent, $\sigma_j$ a dependent set
in $\LL$, both occurring as summands in $d\sigma$, hence $\bigvee
\sigma_i = \bigvee \sigma_j = \bigvee \sigma$. We consider
$\tau=\sigma_i\setminus\{j\}= \sigma_j\setminus\{i\}$; $\tau$ is
independent as a subset of $\sigma_i$, hence it is a maximal
independent set in $\sigma_j$. We see that $\bigvee \tau = \bigvee
\sigma_j = \bigvee \sigma$, in contradiction to $\sigma_i$ being
independent.

\smallskip
\noindent
{\bf (3)} 
{\em $\Psi$ is a quasi-isomorphism.\/} In the case of a geometric
lattice, the classes $[\sigma]$ for independent sets $\sigma$ in $L$
generate $H^*(D,d)$, compare Remark~\ref{rem_Dd}.  Hence, the induced
map $\Psi^*$ is surjective, and, since $H^*(D,d)$ is finite
dimensional, this suffices for $\Psi^*$ to be an isomorphism.
\end{pf}


\section{An outlook}  \label{sect_examples}

With the purpose of going beyond the case of geometric lattices, 
one might be tempted to replace the
notion of independent sets in a geometric lattice with the following
(compare \cite[Sect.\ 3]{Yu99}).

\begin{dfrm} \label{df_independent}
Let $\LL$ be a finite lattice and $\mathfrak A(\LL)$ its set of atoms. 
A subset $\sigma$ in $\mathfrak A(\LL)$ is called {\em independent\/}
if 
\[
       \bigvee \sigma\setminus \{s\} \, \, <\, \,  \bigvee \sigma 
\qquad \mbox{ for any }\, s\in \sigma\, .
\]
\end{dfrm}

We remark again that prior to Section~\ref{sect_formal} we have not been 
referring to any specific property of geometric lattices. 
A careful reading of the proof of Theorem~\ref{thm_Dformal} shows that
there are two points where we had to rely on the intersection lattice 
being geometric.

\begin{quote}
\begin{itemize}
\item[(1)] $H^*(D,d)$ is generated by
  classes $[\sigma]$ of independent sets $\sigma$ in $\LL$.
\item[(2)] The generators
  $\sigma_i$ occurring in $d\sigma$ for $\sigma\,{\subseteq}\,\mathfrak A(\LL)$ 
  are either all independent or all dependent sets in $\LL$.
\end{itemize}
\end{quote}

\begin{explrm} \label{ex_nongeom}
Consider the arrangement~$\Aa$ given by the following four subspaces in $\C^4$:
\[
    U_1\,=\,\{x=u=0\}, \,\,\,
    U_2\,=\,\{y=u=0\}, \,\,\,
    U_3\,=\,\{z=u=0\}, \,\,\,
    U_4\,=\,\{x=y, z=0\}\, .
\]
Its intersection lattice is {\em not\/} geometric, compare the 
Hasse diagram in Figure~{\ref{fig_expllattice}}. 
The independent sets 
according to Definition~\ref{df_independent} are
\[
    1,\ 2,\, 3,\, 4,\, 12,\, 13,\, 14,\, 23,\, 24,\, 34,\, 123\, , 
\] 
where, for brevity, we denote atoms in $\LL$ by their indices.

\begin{figure}[ht]
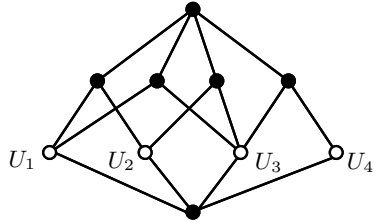

  \begin{picture}(0,0)%
   \includegraphics{expllattice.pstex}%
  \end{picture}%
 \input{expllattice.pstex_t}%

\caption{The intersection lattice of $\Aa$ in Example~\ref{ex_nongeom}}
\label{fig_expllattice}
\end{figure}

Consider the relative atomic complex $(D_{\Aa},d)$.  Its $\LL$-homogeneous
components $D(A)$ consist of a single rank~1 cochain group each, 
for any $A$ in $\LL$ except $A\,{=}\,\hat 0$. The latter reads as follows.
\[
D(\hat 0) \, \, : \quad 0 \, \longrightarrow \,
\langle 1234 \rangle \, \longrightarrow \,
\langle 123, 124, 134, 234 \rangle \, \longrightarrow \,
\langle 14, 24 \rangle\, \longrightarrow \, 0\, ,
\]
with non-trivial cochain groups in degrees $4$, $5$, and $6$. 
Applying $d$, we see that $H^p(D(\hat 0))=0$, unless $p=4$, and 
$H^4(D(\hat 0))=\langle [123] \rangle$.  

Hence, we find that $H^*(D)$ is generated by the classes of the
independent sets in $\LL$. However, contrary to the geometric case,
these classes can be zero, as are $[14]$ and $[24]$ in the present
example.

Since $d(1234)\, = \, 234 - 134 + 124 - 123$ (observe that $123$ is
independent whereas the other classes are dependent!), we have
\begin{equation} \label{eq_explnongeom}
     [123] \, \, = \, \, [234] - [134] + [124]\, .
\end{equation}
This shows that the proof of Theorem~\ref{thm_Dformal} does
not extend to the present arrangement - not all of the classes
on the right hand side of~(\ref{eq_explnongeom}), 
induced by dependent sets in $\LL$, 
can be mapped to $0$ under a quasi-isomorphism between $D$ and
$H^*(D)$.

However, it is easy to check that the arrangement $\Aa$ is formal -
a quasi-isomorphism $D_{\Aa}\,{\rightarrow}\,H^*(D_{\Aa})$ 
can be constructed directly by mapping all $3$ dependent sets above
to the same generating class.

\medskip
A slight variation of this example shows that the proof of 
Theorem~\ref{thm_Dformal} does extend to {\em some\/} non-geometric lattices.
Consider the arrangement of subspaces in $\C^4$ given by
\[
    U_1\,=\,\{x=u=0\}, \,\,\,
    U_2\,=\,\{y=u=0\}, \,\,\,
    U_3\,=\,\{x=y, u=0\}, \,\,\,
    U_4\,=\,\{x=y, z=0\}\, .
\]
Again, $H^*(D)$ is generated by classes of independent sets, and 
$d(1234)\,= \,234-134+124$, hence it is a sum of dependent sets only.
\end{explrm}

\begin{explrm} \label{expl_kequal}
For $n\,{>}\,k\,{\geq}\,2$, consider the 
{\em $k$-equal arrangement\/} $\Aa_{n,k}$  given by the codimension 
$k{-}1$  subspaces in $\C^n$ of the form
\[
   z_{i_1}\, = \,z_{i_2}\, = \,\,  \ldots \,\, = \, z_{i_k}\,, \quad \quad 
   \mbox{ for any }\,\,
   1\leq i_1 < \ldots < i_k \leq n\, . 
\]
Its intersection lattice is the  subposet $\Pi_{n,k}$ of the lattice of set 
partitions of $\{1,\ldots,n\}$ formed by partitions with non-trivial block sizes
larger or equal to~$k$. Observe that $\Pi_{n,k}$ is {\em not\/} geometric for 
$k\,{>}\,2$.

However, in various respects, $k$-equal arrangements do have properties that 
are similar to arrangements with geometric intersection lattices. For example, cohomology 
is generated by classes of independent sets of atoms in~$\Pi_{n,k}$ in the sense of 
Definition~\ref{df_independent}~\cite[Thm.\ 8.8(i)]{Yu02}.  
Our formality proof in Section~\ref{sect_formal}, though, does not extend to 
$D_{\Aa_{n,k}}$; condition (2) mentioned above is violated.
\newline
Consider $\Aa_{7,3}$ and set 
\[
       \sigma = \{ 123, 234, 345, 456, 567\}\, .
\]
Here, the triples $ijk$ are shorthand notation 
for partitions of $\{1,\ldots,7\}$ with only nontrivial block $ijk$. 

The set $\sigma$ is dependent, since removing $234, 345$ or $456$ preserves the
join. We see that $\sigma\setminus \{345\}$ is independent, whereas $\sigma\setminus
\{234\}$ is dependent, e.g., $\bigvee \sigma\setminus\{234,456\}=
\bigvee \sigma\setminus\{234\}$. 
\end{explrm}

%


\end{document}